\documentclass[11pt,twoside]{article}

\usepackage{a4wide,latexsym,longtable,graphicx,subeqn,varioref}
\newcommand{\numsection}[1]{\section{#1}\setcounter{equation}{0}}

\newcommand{\beqn}[1]{\begin{equation}\label{#1}}
\newcommand{\eeqn}{\end{equation}}
\newcommand{\calA}{{\cal A}} \newcommand{\calJ}{{\cal J}}\newcommand{\calC}{{\cal C}}
\newcommand{\eqdef}{\stackrel{\rm def}{=}}
\newcommand{\cvect}[1]{\left( \begin{array}{c} #1 \end{array} \right) }
\newcommand{\tim}[1]{\;\; \mbox{#1} \;\;}
\newcommand{\req}[1]{(\ref{#1})}
\newcommand{\ii}[1]{\{1, \ldots, #1 \}}
\newcounter{algo}[section]
\renewcommand{\thealgo}{\thesection.\arabic{algo}}
\newcommand{\algo}[3]{\refstepcounter{algo}
\begin{center}\begin{figure}[htbp]
\framebox[\textwidth]{
\parbox{0.95\textwidth} {\vspace{\topsep}
{\bf Algorithm \thealgo : #2}\label{#1}\\
\vspace*{-\topsep} \mbox{ }\\
{#3} \vspace{\topsep} }}
\end{figure}\end{center}}
\newcommand{\ms}{\;\;\;\;}

\newcommand{\paperauthor}{Tijana Janji\'{c}\hspace*{-0.02cm}\footnotemark[1],\,\, 
Yvonne Ruckstuhl\footnotemark[1]\,\, and\,
Philippe L.~Toint\footnotemark[2]}
\newcommand{\papertitle}{An algorithm for optimization
   with disjoint linear constraints and its application for predicting rain}

\title{\papertitle}
\author{\paperauthor}
\begin{document}
\renewcommand{\thefootnote}{\fnsymbol{footnote}}
\renewcommand{\thefootnote}{\arabic{footnote}}

\footnotetext[1]{
       Ludwig-Maximilians University Munich,
       Theresienstrasse 37, 80333 M\"{u}nchen, Germany.
       Email: tijana.pfander@lmu.de,}
\footnotetext[2]{
       Namur Centre for Complex Systems,,
       University of Namur,
       61, rue de Bruxelles, B-5000, Namur, Belgium.
       Email: philippe.toint@unamur.be.}

\maketitle

\begin{abstract}
A specialized algorithm for quadratic optimization (QO, or, formerly, QP) with
disjoint linear constraints is presented.  In the considered class of
problems, a subset of variables are subject to linear equality constraints,
while variables in a different subset are constrained to remain in a convex set.
The proposed algorithm exploits the structure by combining steps in the
nullspace of the equality constraint's matrix with projections onto the convex set.
The algorithm is motivated by application in weather forecasting.
Numerical results on a simple model designed for predicting rain show that the
algorithm is an improvement on current practice and that it reduces the  
computational burden compared to a more general interior point QO method.
In particular, if constraints are disjoint and the rank of the set of linear
equality constraints is small, further reduction in computational costs can be
achieved, making it possible to apply this algorithm in high dimensional
weather forecasting problems.
\end{abstract}

\setcounter{page}{1}
\numsection{Introduction}\label{intro-s}

We consider the problem
\beqn{qo-obj}
\min_{x,y} \calJ(x,y) \eqdef (g_x^T, g_y^T)\cvect{x\\y} +
\frac{1}{2} (x^T y^T)\left(\begin{array}{ll} P_{xx} & P_{xy}\\ P_{xy}^T & P_{yy}\end{array}\right)\cvect{x\\y}
\eeqn
subject to
\beqn{qo-constr}
Ax = b \tim{ and } y \leq \ell
\eeqn
where $x$ and $g_x$ belong to $\Re^n$, $y$ and $g_y$ to $\Re^p$, $P_{xx}$ is
an $n \times n$ symmetric real matrix, $P_{yy}$ a $p \times p$ symmetric real matrix,
$P_{xy}$ a $n \times p$ real matrix, $A$ an $m \times n$ real matrix with $m
\leq n$ of rank $m$, $b \in \Re^m$ and the inequality is understood
component-wise. It is easy to extend our discussion to more general bound
constraints where one requires $\ell \leq y \leq u$ for some $\ell, u \in
\Re^p$ with $\ell\leq u$, both of them possibly having infinite components. We
refer to a problem of the type \req{qo-obj}-\req{qo-constr} as having
\emph{disjoint constraints} in the sense that the two sets of constraints of
\req{qo-constr} involve disjoint sets of variables. In what follows, we focus
on the convex case, and assume that
\[
P \eqdef \left(\begin{array}{ll}P_{xx}&P_{xy}\\ P_{xy}^T&P_{yy}\end{array}\right)
\]
is positive definite.

Our motivation for considering problem \req{qo-obj}--\req{qo-constr} arises
from data assimilation in weather forecasting at the convective scale
\cite{Gustetal17}. These forecast are produced by running a complicated
numerical fluid-dynamics model in the future, starting from an initial
condition. This initial condition is however only very partially known and
is typically obtained by fitting available observational data with
short-range model simulation using an optimization algorithm, a technique
known as data assimilation.  For convective scale applications, high
resolution numerical models are used that are very sensitive to the proper
specification of the initial conditions, i.e. predictions of the future can
significantly differ even if the differences in the initial conditions are
small \cite{Lore69}. 

Typical data assimilation algorithms produce the initial condition of the
model by minimizing, perhaps iteratively in a Gauss-Newton framework, an
objective function of type \req{qo-obj}.  Minimization is usually performed
every hour using the new measurements of the atmosphere. The result of the
minimization is a correction to a prediction of the numerical model for a
given time.  In this application, the vector $z=(x,y)$ to be estimated
consists of variables describing the state of the atmosphere at a given time
(such as pressure, temperature, wind direction and speed, \ldots) in all grid
points of the numerical model. Its size $n+p$ ranges between $10^6$ and
$10^9$, resulting in truly large-scale problems. The vector $y$ usually
describes different water phases such as rain, graupel (i.e.\ soft hail) and
snow at all grid points. The physical nature of these variables implies that
they have to be non-negative. In this case, $p$ is approximately one third of
$n$, which remains very large.

In current practice, many data assimilation methods do not
preserve the non-negativity of $y$.  Either these variables are kept constant
if the minimization algorithm attempts to make them negative, or their
optimized values are simply projected onto the positive orthant.  These
techniques clearly interfere with the optimality of the computed result and
the quality of the resulting prediction, especially given the ill-conditioning
of the problem. This was demonstrated by \cite{Janjetal14} to the weather
forecasting community using a simple example where preserving non-negativity
of the $y$ variables (the second part of \req{qo-constr}) as well as the total
mass (the first part of \req{qo-constr}) during the minimization was shown to
be beneficial. Further, \cite{RuckJanj18} extended this observation to
show that similar conclusions hold for a nonlinear multivariable model
designed to test convective scale data assimilation applications.  
These papers use an active set or interior point quadratic optimization (QO) algorithm for solving
the constrained minimization as implemented in matlab \cite{GillMurrWright81,GillMurrSaunWright84} and python \cite{AndeDahlVand10}. Although these results convincingly illustrate
the benefits of including the constraints in the minimization, the QO
algorithm has turned out to be difficult to implement in practice for
applications such as weather forecasting at the convective scale.  This is primarily due to the size of vectors $z$ and
$y$ and the frequency of their estimation that is usually less than an hour.  Both significantly limit the number of iterations of minimization algorithm.  Further, often  for weather forecasting at the convective scale, not only one but rather ensemble of predictions are produced in order to correctly specify, for example, uncertainty of rain at a particular location, even further increasing the computational considerations.  It is the purpose of this short
paper to show that alternative methods do exist and are significantly less
expensive in computer time, thereby making a  practical application affordable.
While these considerations are based on data assimilation for weather
forecasting, we note that the methods discussed here are also applicable to
similar contexts in a wide variety of problems including chemistry, ecosystems
and ocean data assimilation
\cite{BertEvenWack03,SimoBert09,SimoBert12,Buehetal13}, to mention a
few.

In Section~\ref{sec:alg}, we present an algorithm which exploits the fact that
the constraints are disjoint.  Section~\ref{sec:exp} then discusses the results
obtained when applying this algorithm on a known convective-scale example due
to \cite{WuerCrai14}.  In Section~\ref{sec:alg2}, we modify the
algorithm given in Section~\ref{sec:alg} to exploit further properties of our
problem, namely that the equality constraints of \req{qo-constr} are of (very)
low rank, leading to further computational savings.  A discussion and some
perspectives are presented in Section~\ref{sec:con}.

\numsection{The algorithm}\label{sec:alg}

For solving problem \req{qo-obj}--\req{qo-constr}, we propose an active-set
algorithm whose feature is to maintain feasibility with respect to the linear
equality constraints on $x$ at all iterations, while at the same time using classical
projection techniques \cite[Chapter~12]{ConnGoulToin00} to enforce feasibility of the $y$.

If $v$ is a non-negative vector in $\Re^p$ and $\calA \subseteq \ii{p}$, we
denote by $v^\calA$ the vector $v$ reduced to its active component, that
is
\[
v^\calA
= \left\{ \begin{array}{ll} [v]_i & \tim{if } i\in \calA\\
                             \ell &\tim{otherwise,}
          \end{array}\right.
\]
where $[v]_i$ denotes the $i$-th component of $v$.
Similarly, $M^\calA$ is the matrix $M$ reduced to its active columns (and
rows, if it is symmetric).

Our algorithm is stated as Algorithm~\ref{disjqo} \vpageref{disjqo}.

\algo{disjqo}{QO algorithm for disjoint constraints}{
	\begin{description}
		\item[Step 0: Initialization.]
		A feasible starting point $(x_0,y_0)$ is given (i.e., $Ax_0=b$, $y_0\geq
		\ell$), as well as an accuracy threshold $\epsilon >0$. Set $k = 0$.
		\item[Step 1: Active-set update.]
		\beqn{Akdef}
		\calA_k \eqdef \{ i \in \ii{p} \mid [y_k]_i = \ell \tim{and} \nabla_y \calJ(x_k,y_k)
		> 0 \}
		\eeqn 
		\beqn{Akdef}
		\calA_k^c \eqdef \{ i \in \ii{p} \mid i	\not\in \calA_k \}
		\eeqn
		\item[Step 2: Termination test.]
		Terminate if $\| [ \nabla_y \calJ(x_k,y_k) ]^{\calA_k^c} \| \leq \epsilon$.
		\item[Step 3: Search direction computation.]
		Solve
		\beqn{system}
		\left(\begin{array}{ccc} P_{xx} & P_{xy}^{\calA_k^c} & A^T \\*[1ex]
			(P_{xy}^{\calA_k^c})^T & P_{yy}^{\calA_k^c} & 0 \\*[1ex]
			A & 0 & 0 \end{array}\right) \cvect{s_k\\*[1ex]v_k^{\calA_k^c}\\*[1ex]w_k}
		= - \cvect{\nabla_x^{\calA_k^c} \calJ(x_k,y_k)\\*[1ex]\nabla_y^{\calA_k^c}  \calJ(x_k,y_k)\\*[1ex]0}
		\eeqn
		\item[Step 4: Projected search.]
		Determine $\alpha >0 $ such that $(x_{k+1},y_{k+1})$ is the first
		minimizer of 
		$ \calJ \left( x_k+\alpha s_k, \max\Big[y_k+\alpha v_k , \ell \Big]\right)$, where
		$v_k$ is obtained from $v_k^{\calA_k^c}$ by setting $[v_k]_i = 0$ for $i
		\in \calA_k$. 
	\end{description}
}
If not available on the onset, a feasible point can be computed by solving a
linear least-squares problem for $x_0$ and choosing any $y_0 \geq \ell$. 

The third line of \req{system} imposes that $As_k=0$.  It is important that
this equation be satisfied to high precision if exact feasibility with respect
to the linear equality constraint is to be preserved. We refer the reader to
\cite{GoulToin02i} for a discussion of this point. With this caveat, the
system \req{system} can be solved using a Krylov solver like MINRES or GMRES
(see \cite{Saad96} for a description of these methods), or by a ``constrained
preconditioned'' conjugate gradient method (see
\cite{GoulHribNoce01,GoulToin02i}). If this is the case, any preconditioner
must also be reduced (in its $y$ part) to the subset of currently active
variables $\calA_k$.  If dimension and sparsity of $P$ allows (which is
typically not the case in weather forecasting), a stable factorization can also
be used to solve \req{system} accurately.

\numsection{Numerical experiments} \label{sec:exp} 

In order to illustrate the behaviour of the algorithm in our context, we use
the modified shallow water model of \cite{WuerCrai14}. This model has
been used for testing different data assimilation algorithms in
\cite{HaslJanjCrai16, RuckJanj18}. The model
is based on the shallow water (or Saint Venant) equations, which have been
used for a long time in testing both numerical discretization schemes
\cite{Sado75,ArakLamb80, Culletal97,SommNevi09,KeteJaco09} as well as data
assimilation algorithms \cite{CohnParr91,ZengJanj16,Zengetal17}.  As the name
suggests, in \cite{WuerCrai14} the shallow water equations have been 
altered in order to mimic key aspects of convection. To that end, a third
variable rain $r$ was introduced in addition to the velocity (or wind) $u$ and
water height level $h$ fields. The one-dimensional modified shallow water
model consists of following equations:
\beqn{eq:vel} 
\frac{\partial u}{\partial t} + u\frac{\partial u}{\partial x}
+\frac{\partial (\phi + \gamma^2 r)}{\partial x}
= \beta_u + D_u\frac{\partial^2 u}{\partial x^2}, 
\eeqn
with 
\beqn{eq:phi} 
\phi = \left\{ 
\begin{array}{ll} \phi_c  & \tim{if } h > h_c\\
                            gh  &\tim{otherwise,}
\end{array}
\right.
\eeqn
\beqn{modequ2} \frac{\partial r}{\partial t} + u\frac{\partial r}{\partial x}
= D_r\frac{\partial^2 r}{\partial x^2} -\eta r - 
\left\{ 
\begin{array}{ll} \delta \frac{\partial u}{\partial x}, &\tim{if} h>h_r \tim{and} \frac{\partial u}{\partial x} <0 \\
                            0  &\tim{otherwise,}
\end{array}
\right.  \eeqn \beqn{shallowmodel1} \frac{\partial h}{\partial t} +
\frac{\partial (uh)}{\partial x} = D_h \frac{\partial^2 h}{\partial x^2}.
\eeqn
For the physically minded reader, we now briefly describe the meaning of the
various involved quantities. Here, $g$ is the gravity constant and $h_c$
represents the admissible level of free convection. When 
this threshold is reached, the geopotential $\phi$ is set to a lower constant
value $\phi_c$. The parameters $D_u$, $D_r$, $D_h$ are the corresponding
diffusion constants, $\gamma=\sqrt{gh_0}$ is the gravity wave speed for
the absolute fluid layer $h_0$ ($h_0< h_c$).  The small stochastic Gaussian
forcing $\beta_u$ is added at random locations to the velocity, in order to
trigger perturbations and hence convection. Note that this implies that the
location of convection is mostly random. The parameter $\delta$ is the
production rate for rain and $\eta$ is its removal rate. When $h$ reaches the
rain threshold $h_r$ ($h_r > h_c$), rain is 'produced' by adding rain
water mass to the potential, leading to a decrease of the water level and of
buoyancy.

\begin{figure*}
\begin{center}
\includegraphics[width=\textwidth]{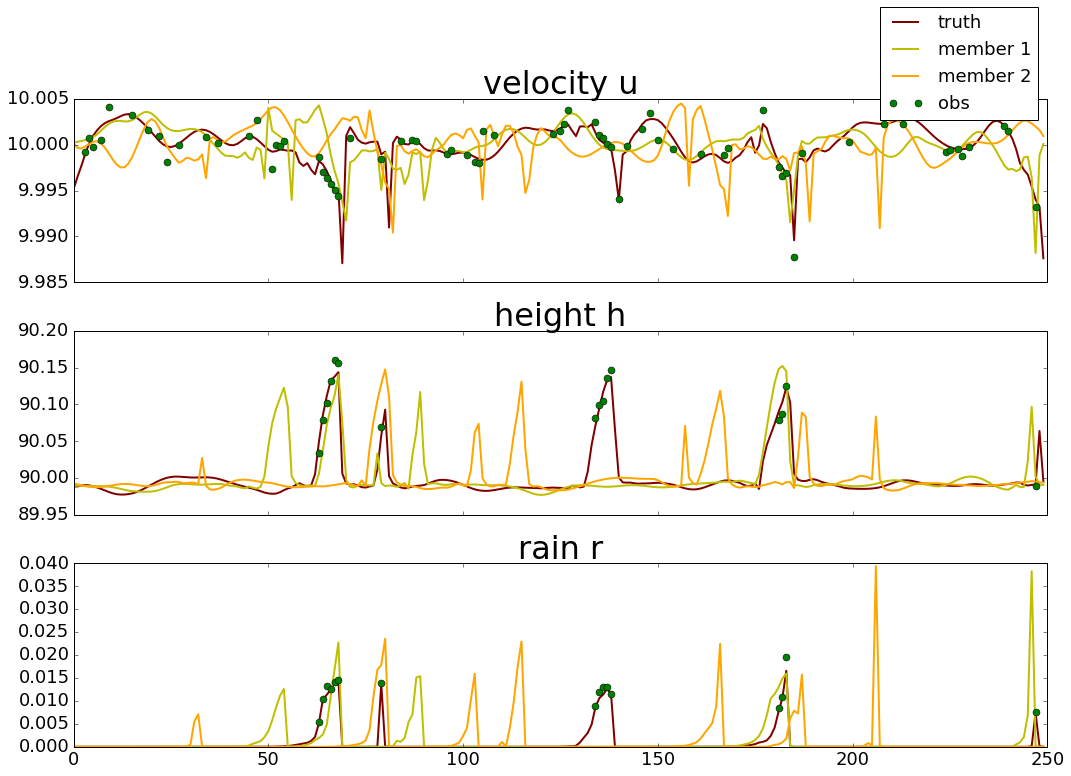}
\caption{Fields produced by running the modified shallow water model with
  three different initial conditions after 60 model time steps.  One of the
  experiments is chosen as nature run (red line).  Observations (green
  circles) are simulated at locations where it rains in the nature run (truth) plus a
  random noise for all fields and, in addition for the $u$ field, an extra
  $25\%$ of   observations are simulated at other locations.  Noise is
  Gaussian for the $h$ and $u$ fields with zero mean  and standard deviations of $0.001$ m/s and
  $0.02$ m, respectively.  Noise is lognormal for the $r$ field, yielding a very
  small observation bias of $0.000825$ and standard deviation of
  $0.00185$.} \label{fig:model}
\end{center}
\end{figure*}

In our numerical implementation of the model, the one dimensional
domain, representing 125 km is discretized with 250 points using standard
second-order centered differences on a staggered grid. To compute evolution
of the dynamical system, the time variable is discretized into time steps of 5
seconds.  The model conserves mass, so the spatial integral over $h$ is
constant in time and the rain $r$ cannot be negative. The other model
parameters are given by
\[
h_0 = 90  \, m,
\ms
h_c = 90.02  \, m,
\ms
h_r = 90.4  \, m,
\ms
D_u = D_h =25000  \, m^{2} s^{-1},
\ms
D_r = 200  \, m^{2} s^{-1},
\]
\[
\phi_c = 899.77 \, m^{2}s^{-2},
\ms
\eta = 2.5 \cdot 10^{-4}  \, s^{-1},
\tim{and}
\delta = 1/300.
\]
The Gaussian stochastic forcing $\beta_u$ has a half width of 4 grid points
and an amplitude of 0.002 \(m/s\). The fields produced by running this model
with three random initial conditions are illustrated in Figure \ref{fig:model}
after 60 model time steps (which is equivalent to five minutes in real
time). As illustrated in Figure \ref{fig:model} position of clouds (height
field) and rain are quite different after only 60 model time steps, mimicking
fast changing convective storms whose intermittency is one of the challenges
of data assimilation on convective scale.  In practice, the availability of
the radar data every 5-15 minutes would preferably be used to recover correct
position and intensities of storms, while longer predictions (48 hours) would
be issued routinely every 6 hours starting from improved initial condition. 

To illustrate Algorithm~\ref{disjqo} we perform a 
twin experiment, where we consider a model run to be the true state
$z=\left({u}^T, {h}^T, {r}^T\right)^T$, which we call the nature run. A vector
of synthetic observations $z^{obs} \in \Re^o$ is then created by randomly
perturbing the nature run such that $z^{obs} =Hz+{\epsilon}$, where $H$ is
the $o \times n$ matrix that determines the location of the observations and
${\epsilon} \in \Re^o$ is a random noise whose components depend on the
observed variable and is computed as follows. For observations of wind and
height field, a Gaussian observation noise is added to the wind $u$ and 
height $h$ fields with zero mean and standard deviations $0.001$ m/s and
$0.02$ m, respectively. A lognormal noise is added to the rain field with
parameters $-8$ and $1.8$, yielding a very small observation bias of
$0.000825$ and standard deviation of $0.00185$. For this choice of parameters,
the observation error for each field is approximately $10\%$ of the maximum
deviation from the field mean. The prior state estimate
$\tilde{z}=({\tilde{u}}^T, {\tilde{h}}^T, {\tilde{r}}^T)^T$ is taken to be
equal to nature run value at a random, much later time.

Given an estimate $\tilde{z}$ and observations $z^{obs}$ of the true state of
the atmosphere, we minimize a quadratic cost function based on the error
covariance matrix of the state estimate $B$ and the observations $R$
respectively, in order to find an improved estimate $z^*=\left({u^*}^T,
{h^*}^T, {r^*}^T\right)^T$ of the true state. We constrain the mass of ${h^*}$
such that ${e}^T{h^*}={e}^T{\tilde{h}}$ and ${r^*}\geq 0$. Specifically the
minimization problem to be solved is  \cite{Lore81,Courtieretal98}:
\beqn{qo-objDA}
\min_{z} \calJ(z)
\eqdef   \frac{1}{2}(z-\tilde{z})^TB^{-1}(z-\tilde{z})
       + \frac{1}{2} (Hz-z^{obs})^TR^{-1}(Hz-z^{obs})
\eeqn
subject to
\beqn{qo-constrDA}
e^Th=e^T\tilde{h} \tim{ and } r \geq 0.
\eeqn
Therefore, in our setup $n=750$, $p=250$, ${x}=(u^T,h^T)^T$ and
$y=r$, $P=B^{-1}+H^TR^{-1}H$  and
$g_x=-B^{-1}\tilde{z}-H^TR^{-1}z^{obs}$. A natural feasible starting point is
$(x_0,y_0)=\left((\tilde{u}^T, \tilde{h}^T)^T, \tilde{r}\right)$.  We estimate
$B$ as a sample covariance from the ensemble of 1000 model simulations that
start from different initial conditions in which correlations that are 10
grids points apart are set to zero.  The observation error covariance matrix
$R$ is taken to be diagonal with values on diagonal corresponding to variances
of distributions used for generating observation error vector ${\epsilon}$. 
We use an LU decomposition with pivoting to solve \req{system} accurately.

\begin{figure*} %[ht!]
\begin{center}
\includegraphics[width=\textwidth]{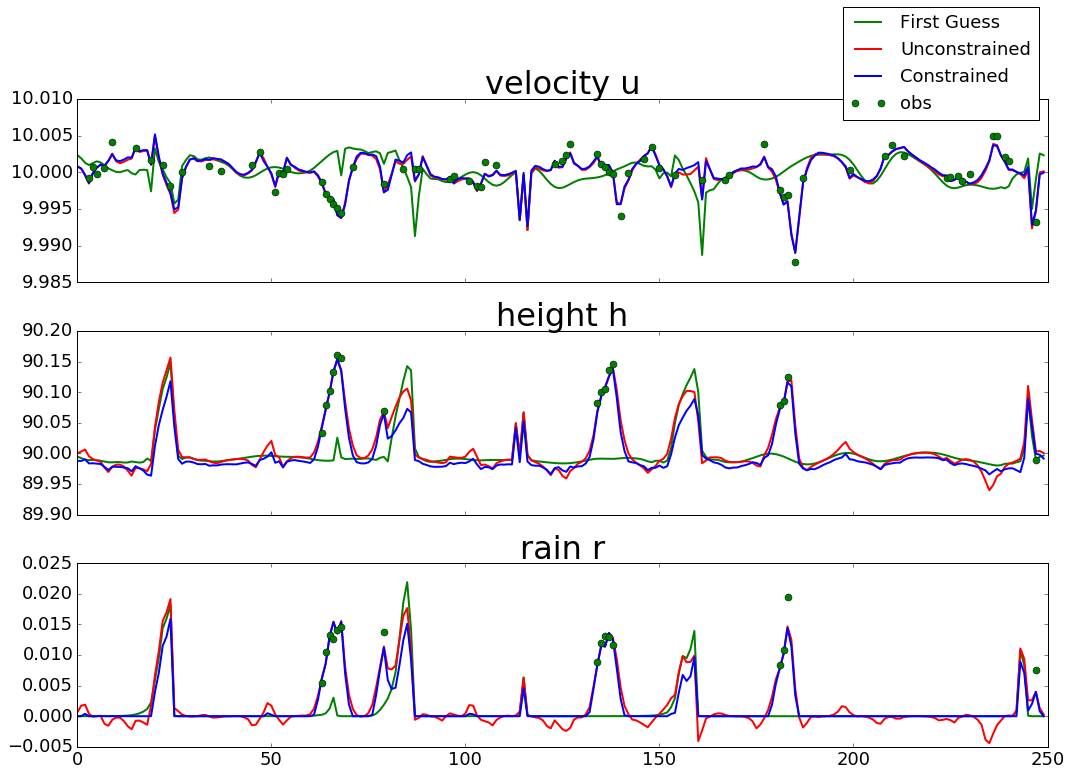}
\caption{Results of minimization for unconstrained (red) and constrained (blue) problems based on prior estimate (green line) and observations (green circles).} \label{fig:estimates}
\end{center}
\end{figure*}

Results for both constrained and unconstrained minimizations are illustrated
in Figure \ref{fig:estimates}.  Constrained minimization produces a slightly
smaller root mean square errors (RMSE) than unconstrained minimization.  In
addition, the value of rain is positive in all grid points.  Although
differences in RMSE between constrained and unconstrained minimization are small after one assimilation cycle, 
in \cite{RuckJanj18} was shown that errors of unconstrained minimization will accumulate over time leading to large
errors in total mass and total rain after repeating data assimilation 250 times, i.e. in less than one
day. Table \ref{tab:par} illustrates the performance of the algorithm.

\begin{table}[h]
\center 
\begin{tabular}{ |p{1cm}||p{3cm}|p{1cm}|p{3cm}|p{1cm}| }
 \hline
 &\hspace*{1cm}$J(z)$     & $|\calA_k^c|$ & $\| [\nabla_{y_i} \calJ(x_k,y_k)]^{\calA_k^c}\|$ & $\alpha_k$ \\
 \hline
 1 & -1.799206e+03 & 133 & 4.811187e+02 & 1.\\
 2 & -1.804615e+03 & 157 & 8.878054e+01 & 1. \\
 3 & -1.805238e+03 & 160 & 1.084530e+01 & 1. \\
 4 & -1.805271e+03 & 162 & 2.373267e-01 & 0.9992  \\
 5 & -1.805271e+03 & 162 & 4.156212e-12 & 1.  \\
 \hline
\end{tabular}
\caption{Illustration of performance of the Algorithm~\ref{disjqo}.}
\label{tab:par}
\end{table}  

As illustrated in Table \ref{tab:par}, Algorithm~\ref{disjqo}  converges
in only five iterations on this example.  If a more general interior point
method like the CVXOPT package \cite{AndeDahlVand10} is applied for
minimization of this problem, the number of iterations required is
typically between ten and twenty.
     
\numsection{The projected algorithm}\label{sec:alg2}

Note that the matrix $A$ in our previous example has a very simple form ($A=(0_u^T,
e_h^T,0_r^T))$ and is of size $1 \times n$.  It obviously has rank one.
We may then easily project the problem into the nullspace of $A$
by defining $Z = I-A^TA/h^2$, the projection onto this nullspace,
and applying the change of variable $x= Z\tilde{x}$ for $\tilde{x} \in
\Re^{n-1}$, which leads to the problem
\beqn{qo-objP}
\min_{\tilde{x},y} \tilde{\calJ}(\tilde{x},y) \eqdef ((Z^Tg_x)^T, g_y^T) \cvect{\tilde{x}\\y} +
\frac{1}{2} (\tilde{x}^T y^T)\left(\begin{array}{ll}
  Z^TP_{\tilde{x}x}Z & Z^TP_{xy}\\
  P_{xy}^TZ         &  P_{yy}
  \end{array}\right)\cvect{\tilde{x}\\y}
\eeqn
subject to
\beqn{qo-constrP}
y \geq \ell.
\eeqn
Problem~\req{qo-objP}-\req{qo-constrP} is now a bound-constrained
quadratic problem, to which standard techniques can be applied, including
for large-scale instances (see \cite{ConnGoulToin92,Dost97,LinMore99}, for example).
A typical such method applies the Conjugate Gradients (CG) \cite{HestStie52} to
minimize the quadratic in the current face, that is the subspace
spanned by the inactive variables at the current iterate, restarting the
procedure as needed when new bound constraints become active during the
calculation and a new (lower dimensional) face must be explored. Note that the
first face contains the negative gradient of the inactive variables, and the
first step of CG performs a (in this case projected) line search along this
direction, yielding what is known as the generalized Cauchy point. Methods differ
essentially by their face changing mechanisms but insist that constraints
active at the Cauchy point are not made inactive during the rest of the
restarted CG steps. A simple version of the resulting algorithm (based on
\cite{ConnGoulToin92}) is now stated as Algorithm~\ref{disjqoP} \vpageref{disjqoP}. 
\algo{disjqoP}{Projected QO algorithm for disjoint constraints}{
\begin{description}
\item[Step 0: Initialization.]
A feasible starting point $(\tilde{x}_0,y_0)$ is given (i.e. $y_0\geq
\ell$), as well as an accuracy threshold $\epsilon >0$.
Compute the projection $Z$ onto the null space of $A$ and set $k = 0$.
\item[Step 1: Active-set update.] Define
\beqn{Akdef}
\calA_k \eqdef \{ i \in \ii{p} \mid [y_k]_i = \ell  \}
\tim{ and }
\calA_k^c \eqdef \ii{p} \setminus \calA_k.
\eeqn
\item[Step 2: Termination test.]
Terminate if the following conditions hold:
\begin{itemize}
\item $\| [\nabla_y \tilde\mathcal{{J}}(\tilde{x}_k,y_k)]^{\calA_k^c} \|\leq \epsilon$
\item $  \| \nabla_{\tilde{x}} \tilde\mathcal{{J}}(\tilde{x}_k,y_k) \|\leq \epsilon$
\item $ \nabla_{y_i} \tilde{\calJ}(\tilde{x}_k,y_k) \geq 0 \tim{ for } i \in	\calA_k$.
\end{itemize}
\item[Step 3: Find the Cauchy point and determine its active set.]
Determine $\alpha >0$ such that $(\tilde{x}_{k}^c,y_{k}^c)$ is the first
minimizer of 
\[
\tilde{\calJ} \left( x_k-\alpha \nabla_{\tilde{x}}
\tilde{\calJ}(\tilde{x}_k,y_k), \max\Big[y_k-\alpha \nabla_{y}
  \tilde{\calJ}(\tilde{x}_k,y_k) , \ell\Big]\right).
\]
Set
\beqn{AkCdef}
\calA_{k,C} \eqdef \{ i \in \ii{p} \mid [y_k^c]_i = \ell  \},
\tim{ and }
\calC_k \eqdef \ii{p} \setminus \calA_{k,C}.
\eeqn
\item[Step 4: Minimization beyond the Cauchy point.]
Apply the CG algorithm to find an approximate minimizer  $(\tilde{x}_{k+1},y_{k+1}^{\calC_k})$ of 
\beqn{qo-objPCG}
((Z^Tg_x)^T, g_{y^{\calC_k}}^T) \cvect{\tilde{x}\\y^{\calC_k}} +
\frac{1}{2} (\tilde{x}^T y^{\calC_k,T})
\left(\begin{array}{ll} Z^TP_{xx}Z        & Z^TP_{xy^{\calC_k}}\\
                        P_{xy^{\calC_k}}^TZ  & P_{y^{\calC_k}y^{\calC_k}}\end{array}\right)
\cvect{\tilde{x}\\y^{\calC_k}}
\eeqn
subject to
\beqn{qo-constrPCG}
y^{\calC_k} \geq \ell^{\calC_k}.
\eeqn
Terminate the CG once one (or more) bound(s) of indices ${j_1,
  \ldots, j_s}$ are encountered, after a maximum number of iterations or once
it has converged. If CG was terminated because bounds were encountered,
restart it after redefining $\calC_k = \calC_k\setminus\{j_1,\ldots,j_s\}$.
Repeat this process until the size of $\calC_k$ does not decrease anymore.
\end{description}
}

In many cases, the efficient application of the CG algorithm requires
preconditioning. We refer the
reader to \cite{GratGuroToin13} for a discussion of suitable strategies in the
context of data assimilation. It is also known that Algorithm~\ref{disjqoP}
could be implemented without Step~4 if mere convergence is wanted, but that
performing conjugate gradient iterations as suggested in
\cite{ConnGoulToin92}, very often significantly reduces number of outer iterations.
This was also observed for our test problem. If subproblems in Step~4 are
solved accurately, Algorithm~\ref{disjqoP} requires three outer iterations, as
illustrated in Table~\ref{tab:par2}.  If the number of conjugate gradient
iterations per outer iteration is fixed a priori (a standard practice in
weather forecasting), the number of outer iterations increases, and could
reach twenty, but the cost of each outer iteration decreases.  The behavior
of the algorithm with the number of CG iterations fixed a priori to $25$,
$50$, $400$ and $800$ is illustrated in Table~\ref{tab:par3}.  For the
computational consideration, we also impose an additional stopping criteria to
Algorithm~\ref{disjqoP}: the algorithm is stopped either when it has converged
or when the number of faces reduces to one. Note that the latter criteria is
met only if the CG iterations for minimising equation (\ref{qo-objPCG}) did
not encounter any bounds prescribed by (\ref{qo-constrPCG}), which would
suggest that the current guess of the active set is fairly accurate, though
not guaranteed to be exact.   Fixing the total number of CG iterations per 
outer iteration limits number of CG restarts during one major iteration and
reduces accuracy as well as cost.  For example, for a fixed number of $25$ CG
iterations, the solutions obtained by Algorithm~\ref{disjqo} and
Algorithm~\ref{disjqoP} only coincide to two 
significant digits, while if $800$ CG iterations are allowed, they share eleven
significant digits.  When allowing the number of CG iterations to increase
from 25 to 800, the total number of iterations performed increases from 300
to 2400 and reaches 2472 in case where no limit is set 
while the cost increases by $75\%$ for no-limit case.  While these number are
encouraging, they also indicate that more attention must be given to
preconditioning. 

\begin{table}[h]
\center 
\begin{tabular}{ |p{0.3cm}||p{2.8cm}|p{0.7cm}|p{3cm}|p{1.6cm}|p{1.2cm}|p{0.7cm}|p{1.7cm}|}
 \hline
 & \hspace*{1cm}$J(z)$ & $|\calA_k^c|$ &
 $\|[\nabla_y\calJ(x_k,y_k)]^{\calA_k^c}\|$
 & \hspace*{0.8cm}$\alpha_k$ & CG its & faces & $\|z_k-z^*\|$ \\
 \hline
 1 &  -1.794983e+03 & 191 & 6.762874e-09 & 6.921e-07 &884 & 58 &  8.521e-02 \\
 2 &  -1.804782e+03 & 167 & 7.144318e-09 & 7.204e-06 & 825 & 9 & 2.193e-02 \\
 3 &  -1.805271e+03 & 162 &  3.945761e-09  & 1.012e-05 & 763 & 6  &  1.572e-12 \\
 \hline
\end{tabular}
\caption{Illustration of performance of the Algorithm~\ref{disjqoP}. In this
  table, ``CG its'' stands for the total number of CG iterations at major
  iteration $k$ and ``faces'' is the number of explored faces at iteration
  $k$. To illustrate the accuracy, the difference is calculated between result
  of each major iteration $z_k$ to $z^*$ an end solution of
  Algorithm~\ref{disjqo} \vpageref{disjqo}.} 
\label{tab:par2}
\end{table} 

\begin{table}[h]
\center 
\begin{tabular}{ |p{0.3cm}||p{2.8cm}|p{0.7cm}|p{3cm}|p{1.6cm}|p{1.2cm}|p{0.7cm}|p{1.7cm}|}
 \hline
 a) & \hspace*{1cm} $J(z)$ & $|\calA_k^c|$ &
 $\|[\nabla_y\calJ(x_k,y_k)]^{\calA_k^c}\|$
 & \hspace*{0.8cm}$\alpha_k$ & CG its & faces & $\|z_k-z^*\|$ \\
 \hline
 1 &    -1.655016e+03  & 140 & 3.574082e+03 & 6.921e-07 &25 & 24  &  2.780e-01 \\
 2 &   -1.791893e+03 & 164 & 1.546087e+03 & 1.411e-06 & 25  & 24  &  2.153e-01\\
 3 &   -1.803668e+03 & 170 &   5.819331e+02 & 8.576e-07 & 25  & 15   &  1.787e-01 \\
 4 &    -1.804859e+03 & 168 &  4.755341e+02 & 9.515e-07 & 25  & 7   &   1.535e-01 \\
 5 &   -1.805179e+03 & 167 & 7.550745e+02 & 7.056e-07 &25  & 7 &  1.033e-01  \\
 6 &    -1.805260e+03 & 164 & 1.648924e+02 & 6.691e-07 & 25 & 3 & 8.738e-02 \\
 7 & -1.805270e+03 & 165 & 1.098210e+02  & 9.226e-07 & 25  & 4 &  5.906e-02  \\
  8 &    -1.755016e+03  & 164 & 1.129395e+02 & 7.539e-07 &25 & 4  &  4.440e-02 \\
 9&   -1.791893e+03 & 163 & 2.762388e+02 & 7.185e-07 & 25  & 3  &  2.991e-02\\
 10 &   -1.803668e+03 & 162 &   3.200715e+01 & 6.718e-07 & 25  & 8   &  2.7e-02 \\
 11 &    -1.804859e+03 & 162 &  5.363492e+01 & 8.069e-06 & 25  & 4   &   1.682e-02 \\
 12 &   -1.805179e+03 & 163 & 1.593856e+01 & 7.072e-07 &25  & 2 &  8.455e-03  \\
 13 &    -1.805260e+03 & 163 & 1.339967e+01 & 6.14e-07 & 25  & 2 & 5.098e-03 \\
 14 & -1.805270e+03 & 163 & 8.526599e+00  & 6.484e-07 & 25 & 2 &  4.336e-03  \\
 15 &    -1.755016e+03  & 163 & 6.887983e+00 & 6.516e-07 &25 & 2  &  3.584e-03 \\
 16 &   -1.791893e+03 & 163 & 3.851198e+00 & 6.514e-07 & 25  & 2  &  2.823e-03\\
 17 &   -1.803668e+03 & 163 &   3.183931e+00 & 7.135e-07 & 25  & 2   &  2.201e-03 \\
 18 &    -1.804859e+03 & 163 &  2.726162e+00 & 7.205e-07 & 25  & 2   &   1.602e-03 \\
 19 &   -1.805179e+03 & 162 & 6.857200e+00 & 9.396e-07 &25  & 1 &  5.692e-04  \\
 \hline
\end{tabular}
\begin{tabular}{ |p{0.3cm}||p{2.8cm}|p{0.7cm}|p{3cm}|p{1.6cm}|p{1.2cm}|p{0.7cm}|p{1.7cm}|}
 \hline
 b) & \hspace*{1cm} $J(z)$ & $|\calA_k^c|$ &
 $\|[\nabla_y\calJ(x_k,y_k)]^{\calA_k^c}\|$
 & \hspace*{0.8cm}$\alpha_k$ & CG its & faces & $\|z_k-z^*\|$ \\
 \hline
 1 &    -1.755016e+03  & 172 & 8.588824e+02 & 6.921e-07 &50 & 41  &  2.292e-01 \\
 2 &   -1.791893e+03 & 172 & 9.320186e+02 & 9.247e-07 & 50  & 22  &  1.571e-01\\
 3 &   -1.803668e+03 & 168 &   1.219848e+02 & 5.457e-07 & 50  & 10   &  6.399e-02 \\
 4 &    -1.804859e+03 & 163 &  2.176264e+02 & 2.483e-06 & 50  & 7   &   3.290e-02 \\
 5 &   -1.805179e+03 & 164 & 5.157843e+01 & 6.450e-07 &50  & 4 &  1.802e-02  \\
 6 &    -1.805260e+03 & 162 & 2.105820e+01 & 6.865e-07 & 50  & 2 & 5.707e-03 \\
 7 & -1.805270e+03 & 162 & 1.431093e+01  & 9.323e-07 & 50  & 1 &  1.158e-03  \\
 \hline
\end{tabular}
\begin{tabular}{ |p{0.3cm}||p{2.8cm}|p{0.7cm}|p{3cm}|p{1.6cm}|p{1.2cm}|p{0.7cm}|p{1.7cm}|}
 \hline
  c) & \hspace*{1cm}$J(z)$ & $|\calA_k^c|$ &
 $\|[\nabla_y\calJ(x_k,y_k)]^{\calA_k^c}\|$
 & \hspace*{0.8cm}$\alpha_k$ & CG its & faces & $\|z_k-z^*\|$ \\
 \hline
 1 &  -1.794983e+03 & 191 & 3.271074e-01 & 6.921e-07 &400 & 58 &  8.522e-02 \\
 2 &    -1.804782e+03 & 167 &  7.592096e-03 & 7.204e-06 & 400 & 9 &2.193e-02  \\
 3 &  -1.805271e+03 & 162 & 2.548443e-03 & 1.012e-05 & 400 & 6  &  2.849e-07 \\
 4 &  -1.805271e+03 & 162 & 2.351854e-08  & 5.207e-07 & 400 & 1  &  4.587e-12 \\
 \hline
\end{tabular}
\begin{tabular}{ |p{0.3cm}||p{2.8cm}|p{0.7cm}|p{3cm}|p{1.6cm}|p{1.2cm}|p{0.7cm}|p{1.7cm}|}
 \hline
 d) & \hspace*{1cm}$J(z)$ & $|\calA_k^c|$ &
 $\|[\nabla_y\calJ(x_k,y_k)]^{\calA_k^c}\|$
 & \hspace*{0.8cm}$\alpha_k$ & CG its & faces & $\|z_k-z^*\|$ \\
 \hline
 1 &  -1.794983e+03 & 191 & 2.929874e-07 & 6.921e-07 & 800 & 58 &  8.521e-02 \\
 2 &   -1.804782e+03 & 167 & 3.247265e-08 & 7.204e-06 & 800 & 9  &  2.193e-02  \\
 3 &  -1.805271e+03 & 162 & 5.351619e-09 & 1.012e-05 & 762 & 6   &  1.583e-12 \\
 \hline
\end{tabular}
\caption{Illustration of performance of the Algorithm~\ref{disjqoP} when the
  maximum number of CG iterations per major iteration is fixed to a) 25, b)
  50, c) 400 and d) 800 respectively.  The notation follows that of Table~\ref{tab:par2}.} 
\label{tab:par3}
\end{table}  

\numsection{Conclusion}\label{sec:con}

We have presented two projection algorithms which exploit the disjoint nature
of constraints typically occurring in weather forecasting applications.
While projection methods may be inefficient when the combinatorial aspect of
selecting the correct active bounds dominate and many faces need to be
explored at each major iteration (in which case the interior-point algorithms
perform better), they do perform well compared to the interior-point
algorithms when the gradient quickly provides a good identification of the
active constraints.  This appears to be the case in our (representative)
application.

The first of our methods, Algorithm~\ref{disjqo}, is more efficient than an
interior point approach on a representative example, but still requires
solving the KKT system \req{system}, which is impractical in weather
forecasting applications due to problem size and frequency of solution.  By
contrast, Algorithm~\ref{disjqoP} exploits the low rank of the linear equality
constraints and uses a well-known iterative approach to compute a possibly
approximate solution while ensuring satisfaction of the constraint. If the
size of the problem is such that the conjugate gradient algorithm is allowed
to converge, the number of outer iterations required by
Algorithm~\ref{disjqoP} is smaller or comparable to that required by
Algorithm~\ref{disjqo}.  If the number of conjugate gradient iterations is
limited from the start (as is often the case in weather forecasting
applications), the number of outer iterations typically increases and finding
the optimal equilibrium between accuracy and cost then depends on the problem
at hand.  A further advantage of Algorithm~\ref{disjqoP} is that its applies
the conjugate gradient to a subproblem whose size is significantly smaller
than that of the KKT system \req{system} (remember that $p\approx n/3$).

The observations made in this note are encouraging (and have already spurred
some interest from the weather forecasting operational centers), but the
authors are aware that adapting the proposed method(s) to a real production
environment remains a significant task, as preconditioning and the details of
the face changing mechanism will need thought and fine tuning.

{\footnotesize
%\bibliographystyle{plain}
%\bibliography{/home/pht/bibs/refs}

\begin{thebibliography}{10}

\bibitem{AndeDahlVand10}
M.~S. Andersen, J.~Dahl, and L.~Vandenberghe.
\newblock Implementation of nonsymmetric interior-point methods for linear
  optimization over sparse matrix cones.
\newblock {\em Mathematical Programming, Series~C}, 2:167--201, 2010.

\bibitem{ArakLamb80}
A. Arakawa and V.~R. Lamb.
\newblock A Potential Enstrophy and Energy conserving scheme for the Shallow Water Equations.
\newblock{\em  Monthly Weather Review}, 109, 18--36, 1980.

\bibitem{BertEvenWack03}
L.~Bertino, G.~Evensen, and H.~Wackernagel.
\newblock Sequential data assimilation techniques in oceanography.
\newblock {\em International Statistical Reviews}, 71:223--241, 2003.

\bibitem{Buehetal13}
M.~B\"{u}hner, A.~Caya, L.~Pogson, T.~Carrieres, and P.~Pestieau.
\newblock A new environment {C}anada regional ice analysis system.
\newblock {\em Atmosphere-Ocean}, 51(1):18--34, 2013.

\bibitem{CohnParr91}
S.~E.  Cohn, and D.~F. Parrish.
\newblock The Behavior of Forecast Error Covariances for a Kalman Filter in Two Dimensions.
\newblock {\em Monthly Weather Review}, 119(8): 1757--1785, 1991.

\bibitem{ConnGoulToin92}
A.~R. Conn, N.~I.~M. Gould, and Ph.~L. Toint.
\newblock {\em {\sf LANCELOT}: a {F}ortran package for large-scale nonlinear
  optimization ({R}elease {A})}.
\newblock Number~17 in Springer Series in Computational Mathematics. Springer
  Verlag, Heidelberg, Berlin, New York, 1992.

\bibitem{ConnGoulToin00}
A.~R. Conn, N.~I.~M. Gould, and Ph.~L. Toint.
\newblock {\em Trust-Region Methods}.
\newblock MPS-SIAM Series on Optimization. SIAM, Philadelphia, USA, 2000.

\bibitem{Courtieretal98}
P. Courtier, E.  Andersson, W.  Heckley, J. Pailleux,  D. Vasiljevic,  M. Hamrud, A. Hollingsworth, F. Rabier, and M. Fisher.
\newblock The ECMWF implementation of three-dimensional variational assimilation (3D-Var). I: Formulation.
\newblock {\em Quarterly Journal of the Royal Meteorological Society}, 127, 1783--1807, 1998.

\bibitem{Culletal97}
M.~J.~P. Cullen, T.  Davies,  M.~H. Mawson, J.~A. James, and S.~C. Coulter.
\newblock An Overview of Numerical Methods for the Next Generation U.K. NWP and Climate Mode1.
\newblock {\em Atmosphere-Ocean}, 35, 425--444,1997.

\bibitem{Dost97}
Z.~Dost\'{a}l.
\newblock Box constrained quadratic programming with proportioning and
  projections.
\newblock {\em SIAM Journal on Optimization}, 7(3):871--887, 1997.

\bibitem{GillMurrWright81}
P.~E. Gill, W. Murray, and M.~H. Wright.
\newblock {\em Practical Optimization}.
\newblock  Academic Press, London, UK, 1981.

\bibitem{GillMurrSaunWright84}
P.~E. Gill, W. Murray, M.~A. Saunders and M.~H. Wright.
\newblock {\em Procedures for Optimization Problems with a Mixture of Bounds and General Linear Constraints}.
\newblock   {\em ACM Trans. Math. Software}, 10:282--298, 1984.

\bibitem{GoulHribNoce01}
N.~I.~M. Gould, M.~E. Hribar, and J.~Nocedal.
\newblock On the solution of equality constrained quadratic problems arising in
  optimization.
\newblock {\em SIAM Journal on Scientific Computing}, 23(4):1375--1394, 2001.

\bibitem{GoulToin02i}
N.~I.~M. Gould and Ph.~L. Toint.
\newblock An iterative working-set method for large-scale non-convex quadratic
  programming.
\newblock {\em Applied Numerical Mathematics}, 43(1--2):109--128, 2002.

\bibitem{GratGuroToin13}
S.~Gratton, S.~G\"{u}rol, and Ph.~L. Toint.
\newblock Preconditioning and globalizing conjugate gradients in dual space for
  quadratically penalized nonlinear-least squares problems.
\newblock {\em Computational Optimization and Applications}, 54(1):1--25, 2013.

\bibitem{Gustetal17}
N.~Gustafsson, T.~Janji\'{c}, C.~Schraff, D.~Leuenberger, M.~Weissman,
  H.~Reich, P.~Brousseau, T.~Montmerle, E.~Wattrelot, A.~Buc\'{a}nek, M.~Mile,
  R.~Hamdi, M.~Lindskog, J.~Barkmeijer, M.~Dahlbom, B.~Macpherson, S.~Ballard,
  G.~Inverarity, J.~Carley, C.~Alexander, D.~Dowell, S.~Liu, Y.~Ikuta, and
  T.~Fujita.
\newblock Survey of data assimilation methods for convective-scale numerical
  weather prediction at operational centres.
\newblock {\em Quarterly Journal of the Royal Meteorological Society},
 144(713):1218--1256, 2018.

\bibitem{HaslJanjCrai16}
M.~Haslehner, T.~Janji\'{c}, and G.~C. Craig.
\newblock Testing particle filters on simple convective-scale models. {Part} 2:
  A modified shallow-water model.
\newblock {\em Quarterly Journal of the Royal Meteorological Society},
  142(697):1628--1646, 2016.

\bibitem{HestStie52}
M.~R. Hestenes and E.~Stiefel.
\newblock Methods of conjugate gradients for solving linear systems.
\newblock {\em Journal of the National Bureau of Standards}, 49:409--436, 1952.

\bibitem{Janjetal14}
T.~Janji\'{c}, D.~McLaughlin, S.~E. Cohn, and M.~Verlaan.
\newblock Conservation of mass and preservation of positivity with
  ensemble-type {K}alman filter algorithms.
\newblock {\em Monthly Weather Review}, 142:755--773, 2014.

\bibitem{KeteJaco09}
G.~S. Ketefian, and M.~Z. Jacobson.
\newblock  A mass, energy, vorticity, and potential enstrophy conserving lateral fluid-land boundary scheme for the shallow water equations.
\newblock {\em Journal of Computational Physics}, 228,1--32, 2009.

\bibitem{LinMore99}
C.~Lin and J.~J. Mor\'{e}.
\newblock {N}ewton's method for large bound-constrained optimization problems.
\newblock {\em SIAM Journal on Optimization}, 9(4):1100--1127, 1999.

\bibitem{Lore69}
E.~N. Lorenz.
\newblock The predictability of a flow which possesses many scales of motion.
\newblock {\em Tellus}, XXI/3:289--307, 1969.

\bibitem{Lore81}
A.~C. Lorenc.
\newblock A global three-dimensional multivariate statistical interpolation scheme.
\newblock {\em Monthly Weather Review}, 109: 701--721, 1981.

\bibitem{RuckJanj18}
Y.~Ruckstuhl and T.~Janji\'{c}.
\newblock Parameter and state estimation with ensemble {K}alman filter based
  approaches for convective scale data assimilation.
\newblock {\em Quarterly Journal of the Royal Meteorological Society},
  144(712): 826--841, 2018.

\bibitem{Saad96}
Y.~Saad.
\newblock {\em Iterative Methods for Sparse Linear Systems}.
\newblock PWS Publishing Company, Boston, USA, 1996.

\bibitem{Sado75}
R.~Sadourny.
\newblock{The dynamics of finite-difference models of the shallow-water equations.}
\newblock {\em Journal of the Atmospheric Sciences}, 1, 119--143,1975.

\bibitem{SimoBert09}
E.~Simon and L.~Bertino.
\newblock Application of the {G}aussian anamorphosis to assimilation in a {3-D}
  coupled physical-ecosystem model of the {N}orth {A}tlantic with the {EnKF}: A
  twin experiment.
\newblock {\em Ocean Science}, 5:495--510, 2009.

\bibitem{SimoBert12}
E.~Simon and L.~Bertino.
\newblock {G}aussian anamorphosis extension of the {DEnKF} for combined state
  and parameter estimation: application to a {1D} ocean ecosystem model.
\newblock {\em Journal of Marine Systems},
  89:1--18, 2012.
  
\bibitem{SommNevi09}
Sommer, M. and N\'evir, P.
\newblock A conservative scheme for the shallow-water system on a staggered geodesic grid based on a Nambu representation.
\newblock {\em Quarterly Journal of the Royal Meteorological Society}, 135, 485--€"494, 2009.

\bibitem{WuerCrai14}
M.~W\"{u}rsch and G.~C Craig.
\newblock A simple dynamical model of cumulus convection for data assimilation
  research.
\newblock {\em Meteorologische Zeitschrift}, 23:483--490, 2014.

\bibitem{ZengJanj16}
Y. Zeng, and T. Janji\'{c}.
\newblock Study of conservation laws with the Local Ensemble Transform Kalman Filter.
\newblock {\em Quarterly Journal of the Royal Meteorological Society},
142 (699): 2359--2372, 2016.

\bibitem{Zengetal17}
Y. Zeng, T. Janji\'{c}, Y. Ruckstuhl, and M. Verlaan.
\newblock Ensemble-type Kalman filter algorithm conserving mass, total energy and enstrophy.
\newblock {\em Quarterly Journal of the Royal Meteorological Society}, 143(708): 2902--2914, 2017.

\end{thebibliography}

}

\end{document}